\documentclass[12pt]{amsart}
\usepackage{amssymb}

\usepackage{mathptmx}

\usepackage{pstcol}
\psset{unit=1.0cm,linewidth=0.4pt}

\def\vertex{\pscircle[fillstyle=solid,fillcolor=black]{0.05}}
\definecolor{light}{gray}{0.9}
\definecolor{medium}{gray}{0.8}

\newcommand{\Z}{\mathbb{Z}}
\newcommand{\R}{\mathbb{R}}
\newcommand{\N}{\mathbb{N}}

\DeclareMathOperator{\pnt}{\raise 0.5mm \hbox{\large\bf.}}

\DeclareMathOperator{\gp}{gp} \DeclareMathOperator{\cn}{cn}
\DeclareMathOperator{\conv}{conv}

\DeclareMathOperator{\relint}{int} \DeclareMathOperator{\Hilb}{Hilb}
 \DeclareMathOperator{\ini}{in}
 \DeclareMathOperator{\rank}{rank}
 \DeclareMathOperator{\Ker}{Ker}

\def\+#1{\relax\ifmmode\if\noexpand #1\relax \mathop{\kern
    0pt^+{#1}}\nolimits\else \kern 0pt^+\!#1 \fi\else$^*$#1\fi}

\newtheorem{thm}{Theorem}
\newtheorem{lem}[thm]{Lemma}
\newtheorem{cor}[thm]{Corollary}

\theoremstyle{definition}

\newtheorem{rem}[thm]{Remark}

\theoremstyle{plain}
\newtheorem*{thm*}{Theorem}

\def\ie{\hbox{i.\,e.}}

\let\phi=\varphi
\title{$h$-vectors of Gorenstein polytopes}

\author{Winfried Bruns}
\address{FB Mathematik/Informatik, Universit\"at Osnabr\"uck, 49069 Osnabr\"uck, Germany}
\email{wbruns@uos.de}

\author{Tim R\"omer}
\address{FB Mathematik/Informatik, Universit\"at Osnabr\"uck, 49069 Osnabr\"uck, Germany}
\email{troemer@uos.de}

\begin{document}

\begin{abstract}
We show that the Ehrhart $h$-vector of an integer Gorenstein
polytope with a regular unimodular triangulation satisfies
McMullen's $g$-theorem; in particular, it is unimodal. This result
generalizes a recent theorem of Athanasiadis (conjectured by
Stanley) for compressed polytopes. It is derived from a more general
theorem on Gorenstein affine normal monoids $M$: one can factor
$K[M]$ ($K$ a field) by  a ``long'' regular sequence in such a way
that the quotient is still a normal affine monoid algebra. This
technique reduces all questions about the Ehrhart $h$-vector of $P$
to the Ehrhart $h$-vector of a Gorenstein polytope $Q$ with exactly
one interior lattice point, provided each lattice point in a
multiple $cP$, $c\in\N$, can be written as the sum of $n$ lattice
points in $P$. (Up to a translation, the polytope $Q$ belongs to the
class of reflexive polytopes considered in connection with mirror
symmetry.) If $P$ has a regular unimodular triangulation, then it
follows readily that the Ehrhart $h$-vector of $P$ coincides with
the combinatorial $h$-vector of the boundary complex of a simplicial
polytope, and the $g$-theorem applies.
\end{abstract}

\maketitle
%
%
%
\section{Introduction}
Let $P \subseteq \R^{n-1}$ be an integral convex polytope and
consider the {\em Ehrhart function} given by $E(P,m)=|\{ z\in
\Z^{n-1} : \frac{z}{m} \in P\}|$ for $m>0$ and $E(P,0)=1$. It is
well-known that $E(P,m)$ is a polynomial in $m$ of degree $\dim(P)$
and the corresponding {\em Ehrhart series} $E_P(t) = \sum_{m\in \N}
E(P,m) t^m$ is a rational function
$$
E_P(t)=\frac{h_0+h_1t+\dots+h_dt^d}{(1-t)^{\dim(P) +1}}.
$$
We call $h(P)=(h_0,\dots,h_d)$ (where $h_d\neq 0$) the
{\em (Ehrhart) $h$-vector} of $P$. This vector was intensively studied in the last
decades (e.g. see \cite{BRHE98} or \cite{ST96}). In particular, the
following questions are of interest:
\begin{enumerate}
\item
For which polytopes is $h(P)$ {\em symmetric}, \ie\ $h_i=h_{d-i}$
for all $i$?
\item
For which polytopes is $h(P)$ {\em unimodal}, \ie\ there exists a
natural number $t$ such that $h_0\leq h_1\leq \dots \leq h_t \geq
h_{t+1} \geq \dots \geq h_d$?
\end{enumerate}

Let us sketch Stanley's approach to Ehrhart functions via
commutative algebra. The results we are referring to can be found in
\cite{BRHE98} or \cite{ST96}. The Ehrhart function of $P$ can be
interpreted as the Hilbert function of an affine monoid algebra
$K[E(P)]$ (with coefficients from an arbitrary field $K$) and where
the monoid $E(P)$ is defined as follows: one considers the cone
$C(P)$ generated by $P\times \{1\}$ in $\R^n$, and sets
$E(P)=C(P)\cap\Z^n$. The monomial in $K[E(P)]$ corresponding to the
lattice point $x$ is denoted by $X^x$ where $X$ represents a family
of $n$ indeterminates. The algebra $K[E(P)]$ is graded in such a way
that the degree of $X^x$ (or of $x$) is the last coordinate of $x$,
and so the Hilbert function of $K[E(P)]$ coincides with the Ehrhart
function of $P$. Since $P$ is integral, $K[E(P)]$ is a finite module
over its subalgebra generated in degree $1$.

However, in general $K[E(P)]$ is not generated by its degree $1$
elements. If it is, then we say that $P$ is \emph{integrally
closed}, and simplify our notation by setting $K[P]=K[E(P)]$.
Evidently $P$ is integrally closed if and only if $E(P)$ is
generated by the integer points in $P\times \{1\}$, or,
equivalently, every integer point in $cP$, $c\in\N$, can be written
as the sum of $c$ integer points in $P$. Our choice of terminology
coincides with that in \cite{BRGUB}. A \emph{unimodular
triangulation} of $P$ is a triangulation into simplices
$\conv(s_0,\dots,s_r)$ such that $s_1-s_0,\dots,s_r-s_0$ generates a
direct summand of $\Z^{n-1}$. If $\dim P\ge 3$, $P$ need not have a
unimodular triangulation. However, if a unimodular triangulation of
$P$ exists, then $P$ is integrally closed; this follows easily from
the fact that a unimodular simplex is integrally closed.

The monoid $E(P)$ is always normal: an element $x$ of the subgroup
of $\Z^n$ generated by $E(P)$ such that $kx\in E(P)$ for some
$k\in\N$, $k\ge 1$, belongs itself to $E(P)$. By a theorem of
Hochster, $K[E(P)]$ is a Cohen--Macaulay algebra. It follows that
$h_i\ge 0$ for all $i=1,\dots,d$. Using Stanley's Hilbert series
characterization of the Gorenstein rings among the Cohen--Macaulay
domains, one sees that $h(P)$ is symmetric if and only if $K[E(P)]$
is a Gorenstein ring. In terms of the monoid $E(P)$, the Gorenstein
property has a simple interpretation: it holds if and only if
$E(P)\cap\relint C(P)$ is of the form $y+E(P)$ for some $y\in E(P)$.
This follows from the description of the canonical module of normal
affine monoid algebras by Danilov and Stanley.

It was conjectured by Stanley that question (ii) has a positive
answer for the {\em Birkhoff polytope} $P$, whose points are the
real doubly stochastic $n\times n$ matrices and for which $E(P)$
encodes the \emph{magic squares}. This long standing conjecture was
recently proved by Athanasiadis \cite{AT}. (That $P$ is integrally
closed and $K[P]$ is Gorenstein in this case is easy to see.)

Questions (i) and (ii) can be asked similarly for the combinatorial
$h$-vector $h(\Delta(Q))$ of the boundary complex $\Delta(Q)$ of a
simplicial polytope $Q$ (derived from the $f$-vector of
$\Delta(Q)$), and both have a positive answer. The Dehn--Sommerville
equations express the symmetry, while unimodality follows from
McMullen's famous $g$-theorem (proved by Stanley \cite{ST80a}): the
vector $(1,h_1-h_0,\dots, h_{\lfloor d/2\rfloor}-h_{\lfloor
d/2\rfloor-1})$ is an $M$-sequence, \ie\ it represents the Hilbert
function of a graded artinian $K$-algebra that is generated by its
degree $1$ elements. In particular, its entries are nonnegative, and
so the $h$-vector is unimodal.

Athanasiadis proved Stanley's conjecture for the Birkhoff polytope
$P$ by showing that there exists a simplicial polytope $P'$ with
$h(\Delta(P'))=h(P)$. More generally, his theorem applies to
compressed polytopes, \ie\ integer polytopes all of whose pulling
triangulations are unimodular. (The Birkhoff polytope is compressed
\cite{ST80, ST96}.) In this note we generalize Athanasiadis' theorem
as follows:

\begin{thm}\label{main_pol}
Let $P$ be an integral polytope such that $P$ has a regular
unimodular triangulation and $K[P]$ is Gorenstein. Then the
$h$-vector of $P$ satisfies the inequalities  $1=h_0 \leq h_1 \leq
\dots \leq h_{\lfloor d/2\rfloor}$. More precisely, the vector
$(1,h_1-h_0,\dots, h_{\lfloor d/2\rfloor}-h_{\lfloor d/2\rfloor-1})$
is an $M$-sequence.
\end{thm}

The regular triangulations are obtained as subdivisions of $P$ into
the domains of linearity of a piecewise affine, concave and
continuous function on $P$, provided this subdivision is really a
triangulation. They can also be defined by weight vectors $(w_x:x\in
P\cap\Z^{n-1})$ in the following way (with the same proviso): one
takes the convex hull $Q$ of the halflines $\{(x,z):x\in
P\cap\Z^{n-1},\ z\ge w_x\}\subseteq\R^n$ and projects the ``bottom''
of $Q$ onto $P$. (See \cite{BRGUB} or \cite{STU96} for a discussion
of regular subdivisions and triangulations.)

Our strategy of proof (whose last step is Lemma \ref{reg_tri} in
Section~3) is to consider the algebra $K[M]$ of a normal affine
monoid $M$ for which $K[M]$ is Gorenstein. (An affine monoid is a
finitely generated submonoid of a group $\Z^n$.) We relate the
Hilbert series of $K[M]$ to that of a simpler affine monoid algebra
$K[N]$ which we get by factoring out a suitable regular sequence of
$K[M]$. In the situation of an algebra $K[P]$ for an integrally
closed polytope $P$, the regular sequence is of degree~$1$, and we
obtain an integrally closed and, up to a translation, reflexive
polytope such that $h(P)=h(Q)$. (However, note that Musta\c{t}\v{a}
and Payne \cite{MP} have given an example of a reflexive polytope
which is not integrally closed and has a nonunimodal $h$-vector.) If
$P$ has even a regular unimodular triangulation, we can find a
simplicial polytope $P'$ such that the combinatorial $h$-vector $h(\Delta(P'))$
of the boundary
complex of $P'$ coincides with $h(P)$. Then it only remains to apply
the $g$-theorem to $P'$.

Without the condition on regularity of the triangulation we can only
find a simplicial sphere $S$ such that $h(P)=h(S)$. If the
$g$-theorem can be generalized from polytopes to simplicial spheres,
then our theorem holds for all polytopes with a unimodular
triangulation.

As a side effect we show that the toric ideal of a Gorenstein
polytope with a square-free initial ideal has also a Gorenstein
square-free initial ideal.

For notions and results related to commutative algebra we refer to
Bruns--Herzog \cite{BRHE98} and Stanley \cite{ST96}. For details on
convex geometry we refer to the books of Bruns and Gubeladze
\cite{BRGUB} (in preparation) and Ziegler \cite{ZI95}.

\section{Gorenstein monoid algebras}
\label{section3}

We fix a field $K$ for the rest of the paper. Let $C$ be a pointed
rational cone in $\R^n$, \ie\ a cone generated by finitely many
integral vectors that does not contain a full line. Such a cone is
the irredundant intersection $\cn(M)=\bigcap_{i=1}^s H_{\sigma_i}^+$
of rational half-spaces. Here $\sigma_i$ is a linear form with
rational coefficients and $H_{\sigma_i}^+=\{x:\sigma_i(x)\ge 0\}$ is
the positive closed halfspace associated with $\sigma_i$. The
hyperplane on which $\sigma_i$ vanishes is denoted by
$H_{\sigma_i}$. Note that for $H_{\sigma_i}$ the form $\sigma_i$ is
unique up to a nonnegative factor. There is a unique multiple with
coprime integral coefficients, and we call this choice of
$\sigma_i$, $i=1,\dots,s$, the {\em support forms} of $C$. The
\emph{extreme integral generators} of $C$ are the shortest integer
vectors in the edges of $C$,

Let $M\subseteq\Z^n$ be a positive affine monoid, \ie\ an affine
monoid whose only invertible element is $0$. Then the cone $\cn(M)$
generated by $M$ is pointed and the map
$$
\sigma: M \to \Z^s,\quad a \mapsto (\sigma_1(a),\dots,\sigma_s(a)),
$$
is injective. It is called the {\em standard embedding} of $M$. It
can be extended to the subgroup $\gp(M)$ of $\Z^n$ generated by $M$,
and we denote the extension also by $\sigma$.

\begin{lem}
\label{cmhelper} Let $M \subseteq \Z^n$ be a positive normal affine
monoid with $\gp(M)=\Z^n$ and $R=K[M]$. Let
$\sigma_1,\dots,\sigma_s$ be the support forms and $\sigma: \Z^n \to
\Z^s$ the standard embedding of $M$. Moreover, let
$\relint(M)=M\cap\relint(\cn(M))$. Then:
\begin{enumerate}
\item
The $\Z^n$-graded canonical module $\omega_R$ is the ideal of $R$
generated
by all $X^z$ for $z \in \relint(M)$.
\item
$R$ is Gorenstein if and only if there exists a (necessarily unique)
$y \in \relint(M)$ such that $\relint(M)=y+M$, and therefore
$\omega_R=(X^{y})$.
\item
$R$ is Gorenstein if and only if there exists a (necessarily unique)
$y \in \relint(M)$ such that $\sigma(y)=(1,\dots,1)$.
\end{enumerate}
\end{lem}
\begin{proof}
(i) and (ii) are well-known results of Stanley and Danilov.
A proof can be found in \cite{BRHE98}.

(iii) Assume that $R$ is Gorenstein. By (ii) there exists $y \in
\relint(M)$ such that $\omega_R=(X^{y})$. We have that $\sigma_i(y)
>0$ for $i=1,\dots,s$ since $y \in \relint(M)$. Fix $i$ and choose
$z \in \relint(M)$ with $\sigma_i(z)=1$. Such an element $z$ can be
found for the following reason. There exists an element $z' \in M$
such that $\sigma_i(z') =0$ and $\sigma_j(z') >0$ for $j \neq i$.
Furthermore there exists $z'' \in \Z^n$ such that $\sigma_i(z'') =1$
by the choice of $\sigma_i$. For $r \gg 0$ the element $z=rz'+z''\in
\relint(M)$ will do the job.

Now $z-y \in M$ and thus $\sigma_i(z-y) \geq 0$. Hence $\sigma_i(y)
\leq 1$ and therefore $\sigma_i(y) =1$. This shows that
$\sigma(y)=(1,\dots,1)$.

Conversely, if there exists $y \in \relint(M)$ such that
$\sigma(y)=(1,\dots,1)$ then it is easy to see that
$\relint(M)=y+M$.

In each case the uniqueness of $y$ follows from the positivity of
$M$.
\end{proof}

Let $M \subseteq \Z^n$ be a positive affine monoid. It is well-known
that $M$ has only finitely many irreducible elements which form the
unique minimal system of generators of $M$. We call the collection
of these elements the {\em Hilbert basis} of $M$, denoted
$\Hilb(M)$. The following is our main result for monoid algebras.

\begin{thm}
\label{mainthm} Let $M \subseteq \Z^n$ be a positive normal affine
monoid and assume that $R=K[M]$ is Gorenstein. Let $y_1,\dots,y_m
\in \Hilb(M)$ such that $\omega_R=(X^{y_1+\dots+y_m})$ is the
$\Z^n$-graded canonical module of $R$. Then:
\begin{enumerate}
\item $X^{y_1}-X^{y_2},\dots,X^{y_{m-1}}-X^{y_m}$ is a regular
sequence for $K[P]$.
\item $S=R/(X^{y_1}-X^{y_2},\dots,X^{y_{m-1}}-X^{y_m})$ is isomorphic to a
Gorenstein normal affine monoid algebra $K[N]$.
\item The canonical module $\omega_{S}$
is generated by the residue class of $X^{y_1}$.
\end{enumerate}
\end{thm}

The reader should note that in general the elements $y_1,\dots,y_m$
are not uniquely determined; in fact, even their number $m$ may not
unique. It is so, however, in the situation of Theorem
\ref{main_pol}.

We isolate the geometric parts of the proof in two lemmas for
which the following notation is useful: for a point $x$ in a
rational cone $C$ with support forms $\sigma_1,\dots,\sigma_s$ we
denote by $\sigma^>(x)$ the set of indices $i$ such that
$\sigma_i(x)>0$. Let $F_i$ denote the facet of $C$ on which
$\sigma_i$ vanishes. Then $i\in\sigma^>(x)$ if and only if $x\notin
F_i$.

Recall that a rational cone $D\subseteq\R^n$ is called unimodular if
it is simplicial (\ie\ spanned by a linearly independent set of
vectors) and its extreme integral generators form a subset of a
basis of $\Z^n$. A unimodular triangulation of a rational cone is a
triangulation into unimodular subcones.

\begin{lem}\label{uni_trian}
Let $C\subseteq\R^n$ be a pointed rational cone with support forms
$\sigma_1,\dots,\sigma_s$ and let $y_1,\dots,y_m\in C\cap \Z^n$ such
that the sets $\sigma^>(y_i)$, $i=1,\dots,m$, form a decomposition
of \{1,\dots,s\} into pairwise disjoint subsets. Furthermore let
$\Gamma$ be the subfan of the face lattice of $C$ consisting of the
faces
$$
F_{j_1,\dots,j_m}=\bigcap_{i=1}^m F_{j_i}, \qquad
j_i\in\sigma^>(y_i),\ i=1,\dots,m,
$$
of $C$ and all their subfaces. Finally, let $\Sigma$ be a
triangulation of $\Gamma$ into rational subcones.
\begin{enumerate}
\item Then
$$
\Delta = \Sigma \cup \bigcup_{j=1}^m \{
\cn(G,y_{i_1},\dots,y_{i_j}): G\in\Sigma,\ 1\leq i_1 <\dots<i_j\leq
m \}
$$
is a triangulation of $C$.
\item If $\sigma(y_1+\dots+y_m)=(1,\dots,1)$ and $\Sigma$ is
unimodular, then $\Delta$ is unimodular.
\end{enumerate}
\end{lem}

\begin{figure}[hbt]
$$
\begin{pspicture}(0,0)(3,2)
\def\Va{0,0}
\def\Vb{3,0}
\def\Vc{0.9,2}
\def\Vd{2.9,1.2}
 \rput(\Va){\vertex}
 \rput(\Vb){\vertex}
 \rput(\Vc){\vertex}
 \rput(\Vd){\vertex}
 \rput(2.95,0.6){\vertex}
 \rput(0.45,1){\vertex}
 \pspolygon[linewidth=1.2pt](\Va)(\Vb)(\Vc)(\Vd)
 \psline(\Va)(\Vc)
 \psline(\Vb)(\Vd)
 \rput(0.1,1){$y_1$}
 \rput(3.3,0.6){$y_2$}
 \rput(2.0,2.0){$\Gamma$}
\end{pspicture}
$$
\caption{The construction of $\Gamma$}\label{Gamma}
\end{figure}
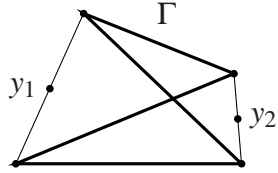

We illustrate the construction of $\Gamma$ by Figure \ref{Gamma} for
the case in which the polytope $P$ is the join $P$ of two line
segments of length $2$ (suitably embedded), $y_1$ and $y_2$ are the
two midpoints (the only possible choice in this case), and
$C=\cn(E(P))$. The bold edges then constitute a ``cross-section'' of
$\Gamma$.

\begin{proof}[Proof of Lemma \ref{uni_trian}]
We may assume that $C$ has dimension $n$. Otherwise we replace
$\R^n$ by $\R C$ and $\Z^n$ by $\R C\cap\Z^n$.

Let us first show that the cones of $\Delta$ are simplicial. It
is enough to consider the maximal elements of $\Delta$. These have
the form $\cn(G,y_1,\dots,y_m)$ for a maximal element $G$ of
$\Sigma$. Let $v_1,\dots,v_r$ be the extreme generators of $G$.
Since $G$ is simplicial, $v_1,\dots,v_r$ are linearly independent.
Assume that
$$
0=\sum_{i=1}^m \lambda_i y_{i} + \sum_{l=1}^r \mu_l v_l\quad
\text{for } \lambda_k, \mu_l \in \R.
$$
Applying $\sigma$ we get $0=\sum_{k=1}^m \lambda_k \sigma(y_{k}) +
\sum_{l=1}^r \mu_l \sigma(v_l)$. By the definition of $\Gamma$ we
have $G\subseteq F_{j_1,\dots,j_m}$ for suitable $j_1,\dots,j_m$.
The hypothesis on $y_1,\dots,y_m$ implies that $\sigma_{j_i}(y_k)=0$
for $k\neq i$, and $\sigma_{j_i}(v_k)=0$ for $k=1,\dots,r$. It
follows that $\lambda_i=0$ for $i=1,\dots,m$, and the linear
independence of $v_1,\dots,v_r$ implies $\mu_l=0$ for $l=1,\dots,r$
as well. Being generated by a linearly independent subset of $\R^n$,
$\cn(G,y_1,\dots,y_m)$ is simplicial.

Next we show that $\Delta$ constitutes a cover of $C$. Let $x \in C$
and set
$$
\lambda_i=\min\biggl\{\frac{\sigma_j(x)}{\sigma_j(y_i)}:
j\in\sigma^>(y_i)\biggr\}\quad \text{for } i=1,\dots,m.
$$
Consider
\begin{equation}
x'=x-\sum_{i=1}^m\lambda_i y_i.\label{multiray}
\end{equation}
First, $\sigma_i(x')\ge 0$ for $i=1,\dots,s$, thus $x'\in C$. (Here
we need the hypothesis on the sets $\sigma^>(y_i)$ in its full
extent!) Second, note that there exists at least one index $j_i$ for
each $i=1,\dots,m$ such that $x\in F_{j_i}$, but $y_i\notin
F_{j_i}$. It follows that $x'$ lies in the face $F_{j_1,\dots,j_m}$
of $C$, and so it belongs to one of the simplicial cones $G$ of
$\Sigma$. Clearly $x\in \cn(G,y_1,\dots,y_m)$.

By definition of $\Delta$, all faces of a cone in $\Delta$ belong to
$\Delta$, too, and it remains to show that the intersection of two
members of $\Delta$ is in $\Delta$. Let $Y_1$ and $Y_2$ be subsets
of $\{y_1,\dots,y_m\}$ and $G_1,G_2$ elements of $\Sigma$. It is
clearly sufficient that
$$
\cn(G_1,Y_1)\cap\cn(G_2,Y_2)=\cn(G_1\cap G_2,Y_1\cap Y_2).
$$
Suppose that $x$ lies in the intersection of the cones. The crucial
point is that both $Y_1$ and $Y_2$ contain the set
$Y'=\{y_i:\lambda_i>0\text{ in }\eqref{multiray}\}$: we have
$$
Y'=\{y_i: \sigma^>(y_i)\subseteq\sigma^>(x)\},
$$
and so $Y'\subseteq Y_1,Y_2$. We conclude that $x'\in
\cn(G_1,Y_1)\cap\cn(G_2,Y_2)$ as well. But since $x'$ belongs to one
of the faces in $\Gamma$, one has $x'\in G_1\cap G_2$ and $x\in
\cn(G_1\cap G_2,Y_1\cap Y_2)$. The converse inclusion is obvious.

It remains to show that the cones in $\Delta$ are unimodular under
the hypothesis of (ii). The unimodularity of $\Sigma$ asserts that
the extreme integral generators of $G$ are part of a basis of $\Z^n$
for every cone $C$ of $\Sigma$. By definition, $G$ is contained in
one of the submodules $L=\bigcap_{i=1}^m H_{\sigma_{j_i}}\cap \Z^n$
of $\Z^n$, which clearly is a direct summand of $\Z^n$ of rank $\ge
n-m$. It is therefore enough that the residue classes of
$y_1,\dots,y_m$ form a basis of $\Z^n/L$ (and, hence, $\rank
L=n-m$). But this is not hard to see : first, the linear forms
$\sigma_{j_i}$ vanish on $L$, and so induce linear forms on
$\Z^n/L$, and, second, the matrix $(\sigma_{j_i}(y_k))$ is the unit
matrix.
\end{proof}

\begin{lem}\label{project}
With the notation and hypothesis of Lemma \ref{uni_trian}(ii) let,
in addition, $V=\R^n/(y_1-y_2,\dots,y_{m-1}-y_m)$, and $\pi:\R^n\to
V$ denote the natural projection. Then the cones $\cn(\pi(G))$ and
$\cn(\pi(G),\pi(y_1))$, $G\in\Sigma$, form a triangulation $\Delta'$
of the cone $\pi(C)\subseteq V$. Moreover, $\Delta'$ is unimodular
(with respect to the lattice $U=\pi(\Z^n)$.
\end{lem}

\begin{proof}
Evidently, the collection $\Delta'$ of the images (i) $\pi(\cn(G))$,
$G\in \Sigma$, and(ii) $\pi(\cn(G,y_1))$, $G\in \Sigma$, covers
$\pi(C)$ since all the images $\pi(y_i)$ coincide with $\pi(y_1)$.
Therefore each of the cones in the triangulation $\Delta$ of $C$ is
mapped onto one of the cones in (i) or (ii).

We have seen in the proof of Lemma \ref{uni_trian} that
$y_1,\dots,y_m$ form part of a basis of $\Z^n$. The same holds for
$y_1-y_2,\dots,y_{m-1}-y_m$. This implies $\Z^n\cap
\Ker\pi=(y_1-y_2,\dots,y_{m-1}-y_m)$ so that $\pi(\Z^n)$ and
$\Z^n/(y_1-y_2,\dots,y_{m-1}-y_m)$ are naturally isomorphic.

Moreover, each cone in $\Delta'$ is generated by part of a basis of
$U$, namely the residue classes of the extreme integral generators
$v_1,\dots,v_r$ of $G$ in case (i) and, in addition, $\pi(y_1)$ in
case (ii). In fact, in the proof of Lemma \ref{uni_trian} we have shown that
$v_1,\dots,v_r,y_1,\dots,y_m$ are part of a basis of $\Z^n$. The
same holds for $v_1,\dots,v_r,y_1-y_2,\dots,y_{m-1}-y_m,y_m$, and we
project onto a submodule generated by the subsystem
$y_1-y_2,\dots,y_{m-1}-y_m$.

The lemma follows once we have proved that $\pi$ maps the union $C'$
of the cones $\cn(G)$ and $\cn(G,y_1)$ bijectively onto $\pi(C)$.
The surjectivity has already been shown above.

Choose $x$ and $x'$ in $C'$ with $\pi(x)=\pi(x')$. Then $x-x'$ is a
linear combination of the differences $y_i-y_{i+1}$. Therefore all
the linear forms $\sigma_{j_1}+\dots+\sigma_{j_m}$ with $j_i\in
\sigma^>(y_i)$, $i=1,\dots,m$ vanish on $x-x'$. There exist
$k_2,\dots,k_m$ and $l_2,\dots,l_m$ with $k_i,l_i\in \sigma^>(y_i)$
and $\sigma_{k_i}(x)=\sigma_{l_i}(x')=0$, $i=2,\dots,m$. Applying
both $\sigma_{j_1}+\sigma_{k_2}+\dots+\sigma_{k_m}$ and
$\sigma_{j_1}+\sigma_{l_2}+\dots+\sigma_{l_m}$ to $x-x'$, one
concludes that $\sigma_{j_1}(x-x')$ must be nonnegative as well as
nonpositive for $j_1\in\sigma^>(y_1)$. So $\sigma_{j_1}(x-x')=0$.
Now the indices in $\sigma^>(y_1)$ are out of the way, and
continuing in the same manner for $y_2,\dots,y_m$ one concludes that
$\sigma_k(x-x')=0$ for all $k$. But then $x-x'=0$, because $C$ is pointed.
\end{proof}

\begin{proof}[Proof of Theorem \ref{mainthm}]
We may assume that $\gp(M)=\Z^n$.
Since $M$ is normal in $\gp(M)=\Z^n$, there exists a unimodular
triangulation of $\cn(M)$ such that each cone in this triangulation
is a unimodular simplicial cone generated by elements of $M$. (See
\cite[Section 2.D]{BRGUB} for a proof of this well-known result.)
Restricting this triangulation to $\Gamma$ we obtain a unimodular
triangulation $\Sigma$ of $\Gamma$.

In view of Lemma \ref{cmhelper}, the cone $C=\cn(M)$ and the elements $y_1,\dots,y_m$ satisfy the
hypothesis of Lemma \ref{uni_trian} (ii). Therefore the triangulation
$\Delta$ of Lemma \ref{uni_trian} (i) is unimodular, and Lemma
\ref{project} provides an induced unimodular triangulation of
$\pi(C)$ where $\pi:\R^n\to V=\R^n/(y_1-y_2,\dots,y_{m-1}-y_m)$ is
the natural projection and the lattice of reference is
$U=\pi(\Z^n)$.

Let $N=\pi(M)$. Since $N$ generates $\pi(C)$ and $C$ has a
unimodular triangulation by elements of $N$, it follows easily that
$N=U\cap \pi(C)$. Therefore $N$ is normal.

As an auxiliary tool we introduce a positive grading on $R$. Let
$k_i=|\sigma^>(y_i)|$, and set
$$
\deg(X^a) = k_2\cdots k_m\sum_{i\in \sigma^>(y_1)}\sigma_i(a) +
\dots + k_1\cdots k_{m-1}\sum_{i\in \sigma^>(y_1)}\sigma_i(a)
$$
for $a \in \Z^n$. The restriction of $\deg$ to $\gp(M)$ is
nonnegative, and for $a\in\gp(M)$ one has $\deg(a)=0$ if and only if
$a=0$. It is obvious that $\deg$ extends to a positive grading of
the ring $R$. Moreover, all preimages of an element of $N$ have the
same degree, so that we have an induced grading on $K[N]$. Finally,
the elements $X^{y_i}-X^{y_{i+1}}$ are all homogeneous of degree
$k_1\cdots k_m$. Therefore, the residue class ring
$$
S=R/(X^{y_1}-X^{y_2},\dots,X^{y_{m-1}}-X^{y_m})
$$
is graded by $\deg$, too.

We want to show that $X^{y_1}-X^{y_2},\dots,X^{y_{m-1}}-X^{y_m}$ is
a regular $R$-sequence. First we prove that $X^{y_1},\dots,X^{y_m}$
is such a sequence. This is most easily seen via the standard
embedding $\sigma:R\to K[\Z_+^s]$ induced by the standard embedding
$\sigma:M\to\Z^s$. Since $K[\Z_+^s]$ is a polynomial ring in $s$
variables and since the sets $\sigma^>(y_i)$ are disjoint, the
monomials $\sigma(X^{y_i})$ ``live'' in pairwise disjoint sets of
variables. So they form a $K[\Z_+^s]$-sequence. Observe that $K[M]$
is a direct summand of $K[\Z_+^s]$ (as a $K[M]$-module) via $\sigma$
(see \cite[Section 4.D]{BRGUB}). Thus $X^{y_1},\dots,X^{y_m}$ is
indeed a regular sequence.

Since
$$
(X^{y_1}-X^{y_2},\dots,X^{y_{m-1}}-X^{y_m},X^{y_m})=(X^{y_1},\dots,X^{y_m})
$$
it is not hard to show that
$X^{y_1}-X^{y_2},\dots,X^{y_{m-1}}-X^{y_m}$ also form a regular
sequence. In fact, all the ideals
$I_{k-1}=(X^{y_1}-X^{y_2},\dots,X^{y_{k-1}}-X^{y_k})$ have height
$k-1$. Since $R$ is Cohen-Macaulay, these ideals are unmixed
(\cite[2.1.6]{BRHE98}), and the next element cannot be a
zero-divisor modulo $I_{k-1}$.

Evidently the ideal $(X^{y_1}-X^{y_2},\dots,X^{y_{m-1}}-X^{y_m})$ is
contained in the kernel of the natural epimorphism $R\to K[N]$, and
in order to show that these two ideals coincide, it is enough that
the rings $S$ and $K[N]$ have the same Hilbert series with respect
to $\deg$.

Recall that the maximal cones in the triangulation $\Delta$ of
$\cn(M)$ and, therefore, all their intersections contain
$y_1,\dots,y_m$. If $C_1,\dots,C_t$ are these maximal cones, then we
see via inclusion--exclusion that
$$
H_R(t) = \sum_{1\leq i \leq t}\  \ \sum_{a \in \gp(M)\cap C_i}
t^{\deg(X^a)} - \sum_{1\leq i<j \leq t}\ \ \sum_{a \in \gp(M)\cap
C_i\cap C_j} t^{\deg(X^a)} \pm \cdots
$$
Let $D_1,\dots,D_t$ be the images of $C_1,\dots,C_t$ with respect to
$\pi$. Then
$$
H_{K[N]}(t) = \sum_{1\leq i \leq t}\ \ \sum_{a \in \gp(N)\cap D_i}
t^{\deg(X^a)} - \sum_{1\leq i<j \leq t} \ \ \sum_{a \in \gp(M)\cap
D_i\cap D_j} t^{\deg(X^a)} \pm \cdots
$$
A comparison of the elements in the unimodular simplicial cones
$C_i$ and $D_i$ yields
$$
H_{K[N]}(t) =(1-t^{k_1\cdots k_m})^{m-1}H_R(t).
$$
But the right hand side of the latter equation is exactly the
Hilbert series of $S$, because
$X^{y_1}-X^{y_2},\dots,X^{y_{m-1}}-X^{y_m}$ is a regular sequence of
$R$, and each element has degree $k_1\cdots k_m$. Hence
$H_S(t)=H_{K[N]}(t)$ and therefore $S \cong K[N]$. Since $S$ is
clearly a Gorenstein ring, its isomorphic copy $K[N]$ is Gorenstein,
too.

It remains to compute the multi-graded canonical module of $S\cong
K[N]$. Since $K[N]$ is Gorenstein, we have to determine the unique
lattice point $q$ in $\relint(N)$ such that $\relint(N)=q+N$ because
then $\omega_{K[N]}=(X^q)$. By construction, $q$ must have degree
$k_1\cdots k_m$, and the residue class of $y_1$ in $U$ is an
interior point of $\cn(N)$ of that degree. This concludes the proof.
\end{proof}

\begin{rem}
The simplicial cone $\cn(y_1,\dots,y_m)$ is the core of the
triangulation $\Delta$ in the sense of \cite{STI} where related
constructions have been discussed. We are grateful to V.~Baty\-rev
for bringing this paper to our attention.
\end{rem}


\section{Gorenstein polytopes}

Let $P \subseteq\R^{n-1}$. We set $E(P,m)=|\{ z\in \Z^{n-1} :
\frac{z}{m} \in P\}|$ and $E(P,0)=1$. In analogy to the rational
function $E_P(t)$ we define
$$
E_{\relint(P)}(t) = \sum_{m\in \N} E(\relint(P),m) t^m
\quad\text{and}\quad E_{\partial(P)}(t) = \sum_{m\in \N}
E(\partial(P),m) t^m.
$$
Observe that $E_{\partial(P)}(t)=E_{P}(t)-E_{\relint(P)}(t)$. In our
situation we have that $E_P(t)=H_{R}(t)$ where $R=K[E(P)]$ and
$E_{\relint(P)}(t)=H_{\omega_{R}}(t)$. Thus these series are
rational with denominator $(1-t)^{\dim(P)+1}$. Moreover,
$E_{\partial(P)}(t)=E_{P}(t)-E_{\relint(P)}(t)=H_{R/\omega_R}(t)$ is
rational with denominator $(1-t)^{\dim(P)}$. (This follows from the
fact that $R/\omega_R$ has Krull dimension equal to $\dim P$.) So it
makes sense to consider the $h$-vectors of these series which we
denote by $h(\relint(P))$ and $h(\partial(P))$. In the following we
present variations and corollaries of Theorem \ref{mainthm}.

\begin{cor}
\label{cor1} Let $P$ be an integrally closed polytope such that
$K[P]$ is Gorenstein. Then there exists a Gorenstein integrally
closed polytope $Q$ such that $\relint(Q)$ contains a unique lattice
point and
$$
h(P) = h(Q) = h(\partial(Q)).
$$
\end{cor}
\begin{proof}
Recall that $R=K[P]$ is the affine monoid ring generated by the
positive normal affine monoid $M=E(P)=C \cap \Z^{n}$ where
$C=\cn((p,1): p \in P)$. Observe that $R$ is $\Z$-graded with
respect to the exponent of the last indeterminate of a monomial and
we will use only this grading for the rest of the proof. All
irreducible elements of $M$ have degree $1$, because $P$ is
integrally closed. Since $R$ is Gorenstein, there exists a unique
lattice point $y \in M$ such that $\relint(M) = y+M$. Choosing
irreducible elements $y_1,\dots,y_m \in M$ such that $y=\sum_{i=1}^m
y_i$ we are in the situation to apply Theorem \ref{mainthm}.

In the proof of the theorem we have constructed the lattice
$U=\gp(M)/(y_i-y_{i+1}:i=1,\dots,m-1)$ and the normal affine lattice
monoid $N \subseteq \gp(M)$ such that $K[N]$ is Gorenstein. The
monoid $N$ is also homogeneous with respect to the grading induced
by that of $M$ and generated by the degree 1 elements. Thus it is
polytopal by \cite[Proposition 1.1.3]{BrGT}, and $K[N]=K[Q]$ for the
polytope $Q$ spanned by the degree 1 elements of $N$. It has also
been shown that the canonical module of $K[Q]$ is generated by a
degree 1 element, the residue class of $X^{y_1}$, which we denote by
$X^p$. Thus $Q$ can have only one interior lattice point, namely
$p$. The $h$-polynomial of $K[P]$ and the one of $K[Q]$ coincide
since $K[Q]\cong K[P]/(X^{y_i}-X^{y_{i+1}}, i=1,\dots,m-1)$ and
$X^{y_1}-X^{y_{2}},\dots, X^{y_{m-1}}-X^{y_{m}}$ is a regular
sequence homogeneous of degree $1$.

It follows from
$$
E_{\partial(Q)}(t) = E_{Q}(t)-E_{\relint(Q)}(t) =
H_{K[Q]}(t)-H_{\omega_{K[Q]}}(t) = H_{K[Q]}(t)-t\cdot H_{K[Q]}(t)
$$
that $h(Q) = h(\partial(Q))$. For the last equality we have used the
fact that $\omega_{K[Q]}=(X^p) \cong K[Q](-1)$ with respect to the
considered grading. This concludes the proof.
\end{proof}

Up to a translation the Gorenstein polytopes with an interior
lattice point are exactly the \emph{reflexive polytopes} used by
Batyrev in the theory of mirror symmetry; see \cite{BA94}. Therefore
the previous corollary reduces all questions about the $h$-vector of
integrally closed Gorenstein polytopes to integrally closed
reflexive polytopes. However, as shown by Musta\c{t}\v{a} and Payne
\cite{MP}, there exist reflexive polytopes that are not integrally
closed and whose $h$-vector is not unimodal.

If $S$ is a simplicial sphere (or even the boundary of a simplicial
polytope), then we can speak of its combinatorial $h$-vector (which
one can read as the $h$-vector of the Ehrhart series of the
geometric realization of $S$ in the boundary of a suitable unit
simplex.)

\begin{cor}
\label{cor2} Let $P$ be an integer polytope such that $K[P]$ is
Gorenstein.
\begin{enumerate}
\item If $P$ has a unimodular triangulation, then there exists
a simplicial sphere $S$ such that $h(P)=h(S)$.
\item If $P$ has a regular unimodular triangulation, then there
exists a simplicial polytope $P'$ such that $h(P)=h(\Delta(P'))$
\end{enumerate}
\end{cor}
\begin{proof}
Polytopes with a unimodular triangulation are integrally closed. So
we can proceed as in the proof of Corollary \ref{cor1} and use the
same notation. The only change is that we start with the given
(regular) unimodular triangulation $\Xi$ of $P$. It induces a
unimodular triangulation of $\cn((p,1): p\in P)$ from which we
derive the triangulation $\Delta$ of $\cn((p,1): p\in P)$ as in the
proof of Theorem \ref{mainthm}. Thus the simplicial cones in
$\Delta$ have generators of degree~$1$, and so it induces a
unimodular triangulation $\Delta_1$ of $P$. The restriction of
$\Gamma$ to $P$ is denoted by $\Gamma_1$: the faces in $\Gamma_1$
are the intersections of the faces of $\Gamma$ with $P$. Then
$\Gamma_1$ is a subcomplex of $\partial P$, and $\Xi|\Gamma_1$ is a
subcomplex of $\Delta_1$. More precisely,
$\Delta_1=(\Xi|\Gamma_1)*\delta$ where $\delta$ is the simplex
generated by the lattice points in $P$ representing the irreducible
elements $y_1,\dots,y_m$ in Theorem \ref{mainthm}. For simplicity we
denote the lattice points also by $y_1,\dots,y_m$.

With the notation of the proof of Theorem \ref{mainthm}, $\Delta$
induces a unimodular triangulation $\Delta'$ of $\cn(N)$ with
generators of degree~$1$ and thus a unimodular triangulation
$\Delta'_1 $ of the (integrally closed) integer polytope $Q$.

Moreover, $K[Q]$ is Gorenstein, $h(P)=h(Q)=h(\partial(Q))$ and
$\relint(Q)$ contains a unique lattice point $p$. For (i) we simply
choose $S=\partial Q$ with triangulation $\Delta'_1|\partial Q$.

For (ii) we first show that the triangulation $\Delta_1$ is regular
since this fact will be needed for an application to initial ideals.
For the same reason we use Sturmfels' correspondence between
monomial initial ideals and regular triangulations of $P$
\cite[Ch.~8]{STU96}. Since $\Xi$ is a regular unimodular
triangulation of $P$, there exists a weight vector $w=(w_x:x\in
P\cap\Z^{n-1})$  such that (i) $\Xi$ is the regular subdivision of
$P$ induced by $w$, and (ii) the initial ideal $\ini_{w}(I_P)$ of
the toric ideal $I_P$ is the Stanley-Reisner ideal of $\Xi$ (as an
abstract simplicial complex). By adding constants to the weights and
scaling them simultaneously we can assume
$$
1\le w_x <1+ \frac 1{n} \quad \text{for all }x\in P\cap\Z^{n-1}.
$$
The toric ideal $I_P$ lives in the polynomial ring $T=K[Y_x:x\in
P\cap\Z^{n-1}]$. It is just the kernel of the natural epimorphism
$\phi:T\to K[P]$, sending $Y_x$ to the monomial $X^{(x,1)}$. It is
generated by binomials that are homogeneous with respect to the
standard grading on $T$.

We define a new weight vector $w'$ on $P\cap \Z^{n-1}$ by keeping
the weight $w_x$ for $x\notin \{y_1,\dots,y_m\}$ and setting $w'_y=0$
for $y_1,\dots,y_m$.

Let $J$ be the Stanley-Reisner ideal of $\Delta_1$. It is enough to
show that $J$ is contained in the initial ideal $\ini_{w'}(I_P)$:
first, $J\subseteq\ini_{w'}(I_P)$ implies $J=\ini_{w'}(I_P)$ since
both residue class rings $T/J$ and $T/\ini_{w'}(I_P)$ have the same
Hilbert function, namely $E(P,{-})$ (we use the unimodularity of
$\Delta_1$). Second, the equality $J=\ini_{w'}(I_P)$ implies that
the regular subdivision of $P$ induced by the weight vector $w'$ is
exactly $\Delta_1$.

Let $M\subseteq P\cap\Z^{n-1}$ and denote the product of the
indeterminates $Y_x$, $x\in M$, by $Y^M$. Then the ideal $J$ is
generated by all monomials $Y^M$ such that $\conv(M)$ is not a face
of $\Delta_1$ and is minimal with respect to this property. In
particular, $Y^M$ has at most $n+1$ factors (for reasons of
dimension). The crucial point is that no such $M$ can contain a
point $y\in P\setminus |\Gamma_1|$, as follows from the construction
of $\Delta_1$.

Suppose first that $\delta=\conv(M)$ is contained in one of the
faces belonging to $\Gamma_1$. Then  $\delta$ is a nonface of $\Xi$
and therefore a nonface of $\Delta_1$ since both triangulations
agree on $\Gamma_1$.

Suppose second that $\delta$ is not contained in a face of
$\Gamma_1$. Then the barycenter of $\delta$ lies in the interior of
one of the unimodular (!) simplices in $\Delta_1$ that have one of
the lattice points $y_i$ as a vertex. Therefore the epimorphism
$\phi:T\to K[P]$ maps $Y^M$ to a monomial that can also be
represented as a monomial involving one of the variables $Y_{y_i}$.
However, this second monomial has the same total degree and strictly
smaller weight with respect to $w'$, as the reader may check. So
$Y^M$ appears in the initial ideal.

This concludes the proof of the regularity of $\Delta_1$. That
$\Delta_1'$ is regular, is seen in the same way. The only difference
is that the lattice points $y_1,\dots,y_m$ are identified to a
single one. For (ii) it remains to apply the next lemma.
\end{proof}

We include a lemma on regular triangulations that we have not found
in the literature.

\begin{lem}\label{reg_tri}
Let $Q\subseteq\R^{n-1}$ be a polytope with a regular triangulation
$\Sigma$. Then there exists a simplicial polytope $P'$ such that the
boundary complex of $P'$ is combinatorially equivalent to
$\Sigma|\partial P$.
\end{lem}

\begin{proof}
We choose a convex, piecewise affine function $f:Q\to\R$ such that
$\Sigma$ is the subdivision of $Q$ into the domains of linearity of
$f$. We can assume that $f(x)>0$ for all $x\in Q$. Consider the
graph $G$ of $f$ in $\R^n$. Then $G$ is a polytopal complex whose
faces project onto the faces of $\Sigma$.

We choose a point $(x,z)\in\R^n$, $x\in\relint(Q)$ and $z\ll 0$ such
that $(x,z)$ lies ``below'' all the hyperplanes through the facets
of $G$. Then we form the set $C$ as the union of all rays emanating
from $(x,z)$ and going through the points of $G$. It is not hard to
check that $C$ is in fact convex: let $a,b\in C$ and consider a
point $c$ on the line segment $[a,b]$; we have to show that the ray
from $(x,z)$ through $c$ meets $G$. We may assume that $a,b\in G$.
Let $c'=(c_1,\dots,c_n)$ be the projection of $c$ in $\R^n$ along
the vertical axis. By convexity of $f$ one has $f(c')\le c_{n+1}$,
and, for the same reason, the graph of $f$ over the line segment
$[x,c']\subseteq Q$ lies below the line segment
$\bigl[(x,f(x)),(c',f(c')\bigr]$. It follows that the line segment
$\bigl[(x,z),c]$ intersects the graph of $f$.

The decomposition of $\partial G$ (as a manifold with boundary) into
maximal polytopal subsets is combinatorially equivalent to the
collection of the maximal simplices in $\Sigma|\partial Q$. On the
other hand, it is also combinatorially equivalent to the collection
of the facets of the cone $C$ (with apex in $(x,z)$). (This requires
an argument very similar to the one by which we have proved the
convexity.) Therefore we obtain the desired polytope $P'$ as a
cross-section of $C$.
\end{proof}

With Corollary \ref{cor2}(ii) the proof of Theorem \ref{main_pol} is
complete since the $g$-theorem applies to $h(\Delta(P'))$.

We are grateful to Ch.\ Haase for pointing out to us that the
hypothesis of regularity cannot be omitted in Corollary
\ref{cor2}(ii) and for suggesting the proof of Lemma \ref{reg_tri}.
The assumptions of Corollary \ref{cor2} appear at several places in
algebraic combinatorics as has been discussed in \cite{AT}.

\begin{rem}
In some special situations we can omit the assumption that $P$ has a
(regular) unimodular triangulation and obtain directly from
Corollary \ref{cor1} that the $h$-vector of $P$ is unimodal. More
precisely, assume that $\dim(P) \leq m+4$. Then $\dim(Q) \leq 5$ and
it follows from a result of Hibi \cite{HI} that the $h$-vector of
$Q$ is unimodal.
\end{rem}

We conclude by drawing a consequence for the toric ideal $I_P$ of
$P$. As we have seen in the proof of Corollary \ref{cor2} its
initial ideal with respect to the weight vector $w'$ is the
Stanley-Reisner ideal $J$ of the simplicial complex
$\Delta_1*\delta$ (notation as in the proof of Corollary
\ref{cor2}). Since $\Delta_1$ is combinatorially equivalent to the
boundary of a simplicial polytope, it follows that $R/J$ is a
Gorenstein ring, and we obtain:

\begin{cor}\label{in_id}
Let $P$ be an integer Gorenstein polytope such that the toric ideal
$I_P$ has a squarefree initial ideal. Then it also has a square-free
initial ideal that is the Stanley-Reisner ideal of the join of a
boundary of a simplicial polytope and a simplex, and thus defines a
Gorenstein ring.
\end{cor}

The corollary answers a question of Conca and Welker, and the
methods of this note were originally designed for its solution.
See \cite[Question 6]{COOB} and \cite{COPRE} for more details
related to this result.

\end{document}